\documentclass{article}
\usepackage{amssymb,amsmath,graphicx,ucs,hyperref}
\usepackage[utf8x]{inputenc}

\newtheorem{question}{Open question}
\newtheorem{conjecture}{Conjecture}

\title{Packing Unequal Disks in the Euclidean Plane}
\author{Thomas Fernique}

\date{}

\begin{document}
\maketitle

\begin{abstract}
A packing of disks in the plane is a set of disks with disjoint interiors.
This paper is a survey of some open questions about such packings.
It is organized into five themes: compacity, conjugacy, density, uniformity and computability.
\end{abstract}

\section{Compacity}
\label{sec:1}

A packing of disks is said to be {\em compact} or {\em triangulated} if its {\em contact graph}, that is, the graph which connects the center of tangent disks, is a triangulation.

When the disks are all of the same size, there is a unique compact packing, known as the {\em hexagonal compact packing} (HCP): the disks are centered on a triangular grid and each disk is surrounded by $6$ tangent disks.

\begin{figure}[hbt]
\centering
\includegraphics[width=0.8\textwidth]{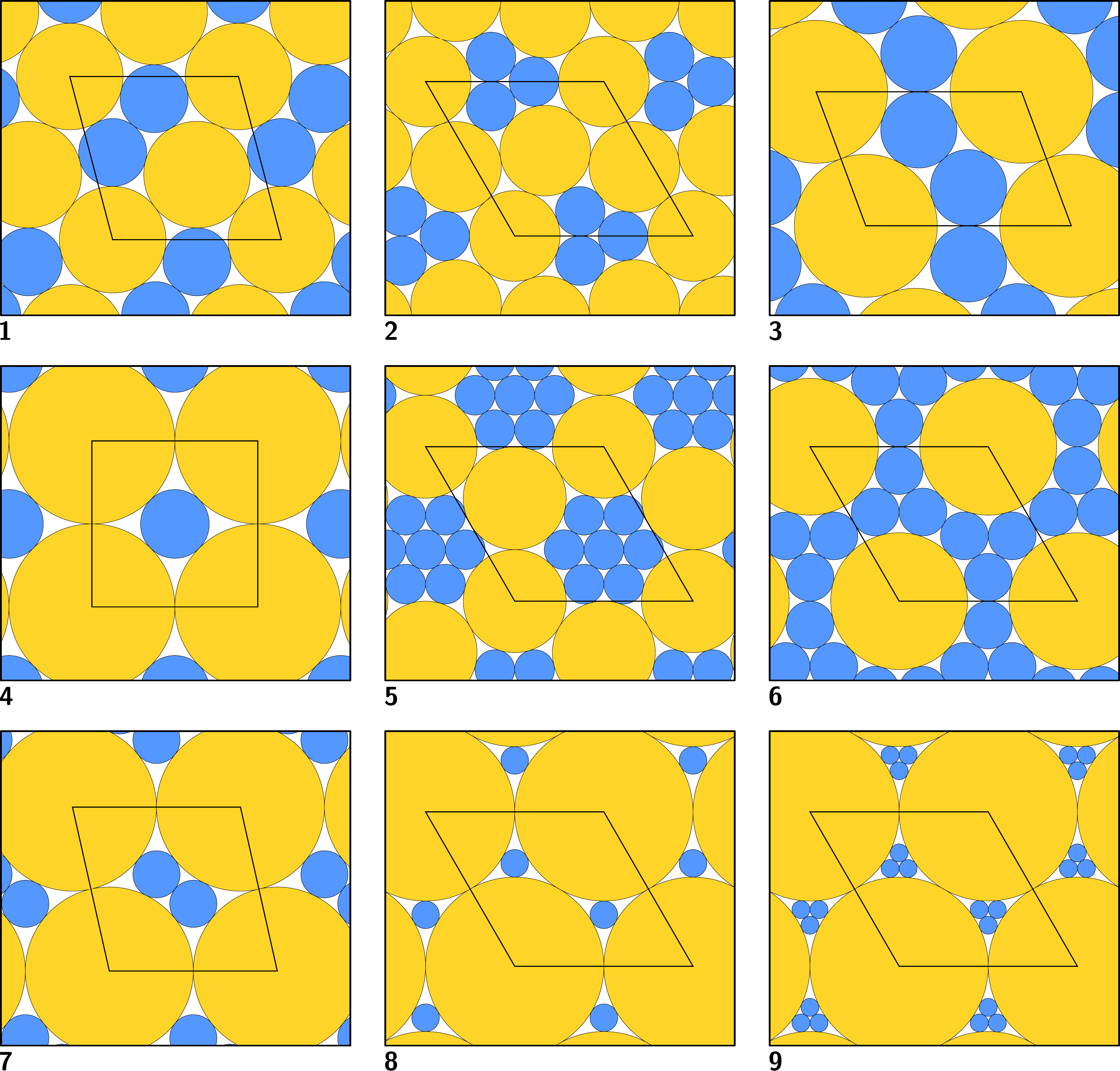}
\caption{
An example of a triangulated binary packing for each of the $9$ possible ratios.
Each packing is periodic and the black parallelogram depicts a fundamental domain.
}
\label{fig:1}
\end{figure}

The situation becomes more interesting for packings with two sizes of disks, called {\em binary packing}.
For example, with a well chosen ratio, it is possible to insert a small disk between three disks of a HCP so that it is tangent to these three disks.
There are several other folk constructions (see, e.g., \cite{LH93}), but it is only in 2006 that the ratios that allow a triangulated packing were characterized by Kennedy \cite{Ken06}: there are $9$ in all (Fig.~\ref{fig:1}).

The proof of \cite{Ken06} relies on the fact that for a triangulated packing to be possible, the sizes of the disks must allow to form a {\em corona}, i.e. a sequence of disks each tangent to the next and all tangent to a central disk.
The sum of the angles of the triangles of the contact graph of such a corona having to sum to $2\pi$, taking the cosine and performing some trigonometric manipulations gives an integer polynomial whose roots define the disk sizes.
The proof is rather simple because the combinatorial of possible coronas is quite limited and the corresponding (univariate) polynomials are of small degree (at most $8$ for the case $2$ in Fig.~\ref{fig:1}).
The same approach extends to more disk sizes but becomes much more complex:
\begin{itemize}
\item finding the coronas - defined by the sequence of sizes of their disks - becomes challenging as their number grows rapidly (in $k^n$ where $k$ is the number of different sizes and n is the number of disks in the corona - which is not even clear that it is bounded for $k$ fixed);
\item finding the sizes - defined as roots of a system of (multivariate) polynomials associated with coronas - is not strictly necessary, but at least we must be able to determine the set of all coronas compatible with these sizes;
\item finding a packing - assuming that all the coronas are known - could be an intractable algorithmic problem (we shall come back to this point in more detail in Section~\ref{sec:5}).
\end{itemize}
This leads to the following challenge:

\begin{question}
How many values $0 < r_1 < \ldots < r_n = 1$ allow a triangulated packing by disks of size $r_1 ,\ldots , r_n$ where each size appears at least once?
\end{question}

In \cite{Mes20}, Messerschmidt was able to bound by $13617$ the number of pairs $(r, s)$ such that there exists a triangulated packing of disks of size $1 > r > s$.
The exact number turned out to be $164$ \cite{FHS21} - establishing this result required several tricks to get around the problems mentioned above.
Few things are known for more disks, with the notable exception that there is always a finite number of such tuples \cite{Mes23}, which is not obvious because the system of equations associated with coronas might not be zero-dimensional.

We could also discard the packings which are derived from packings with fewer disk sizes by inserting disks into holes: they are less interesting and often require more computation (because of the small size of the disks and thus the large size of the coronas).
This would, for example, eliminate the last two cases in Fig.~\ref{fig:1} and $34$ from the $164$ ternary cases.

\section{Conjugacy}
\label{sec:2}

In \cite{Ken06}, Kennedy does not only characterize all the ratios which allow a triangulated packing but also describes the set of possible packings for each ratio (except the one numbered $2$ in Fig.~\ref{fig:1}).
For such a description, it is convenient to introduce some terminology.

First, let us introduce some notions from {\em tiling theory}.
A {\em tile} is a polygon with a color map defined on each of its points (actually only the boundary is relevant).
Tiles can be duplicated, translated or rotated.
Two tiles are said to {\em match}  if their interior do not overlap and whenever a point belongs to both boundaries, the corresponding color maps give it the same color.
Then, a {\em tiling} is a covering of the Euclidean plane by matching tiles.
In other words, a tiling is nothing more than an infinite jigsaw puzzle!

\begin{figure}[hbt]
\centering
\includegraphics[width=0.8\textwidth]{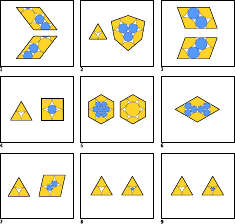}
\caption{
The tilings defined by each of these tile sets correspond to the sets of all the possible triangulated binary packings (Fig.~\ref{fig:1}).
Can you visualize them?
}
\label{fig:2}
\end{figure}

Every triangulated packing by given disk sizes can be seen as a tiling: the tiles are the triangles of the contact graph ($n$ disk sizes yield $n(n^2 + 2)/3$ tiles, or $n(n + 1)(n + 2)/6$ up to reflection).
However, we can usually give a smaller set of tiles which often makes it easier to get an intuition of the whole set of packings.
Fig.~\ref{fig:2} illustrates this.

Further, let us introduce some notions from {\em symbolic dynamics}.
A {\em tiling space} is the set of tilings which can be formed by a given set of tiles.
We say that the tiling space defined by a tile set $\mathcal{T}$ {\em factors} onto the tiling space defined by a tile set $\mathcal{T}'$ if there is a map $\phi : \mathcal{T}\to\mathcal{T}'$ such that, for any tiling by tiles in $\mathcal{T}$ , replacing each tile by its image by $\phi$ yields a tiling (that is, the new tiles match).
If the latter tiling space also factors onto the former, then both are said to be {\em conjugated}.
Intuitively, they are identical up to some ``flourish''.

For example, in Fig.~\ref{fig:2}, cases $8$ and $9$ are conjugated: map the tile without small disk onto itself and the tile with one small disk onto the tile with three small disks (or conversely).
Both factors on HCP (the hexagonal compact packing): map the two tiles on the tile without small disk.
Case $5$ also similarly factors onto HCP, while case $6$ is conjugated to HCP (there is only one possible tiling).
All these cases ($5$, $6$, $8$, $9$ and HCP) are called {\em hexagonal} in \cite{FHS21}.
What about the other cases?

\begin{figure}[hbt]
\centering
\includegraphics[width=0.8\textwidth]{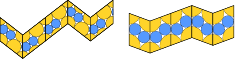}
\caption{Cases $1$ and $3$ of Fig.~\ref{fig:2}: the two tiles freely alternate and form an infinite horizontal stripe which can then be repeated in the vertical direction.}
\label{fig:3}
\end{figure}

Both cases $1$ and $3$ are similar (Fig.~\ref{fig:3}).
Such cases are called {\em laminated} in \cite{FHS21}.
They are not conjugated because the tiles have different shapes, but they do if we slightly extend the notion of factor as follows: after applying $\phi$, tiles can be moved so that two tiles are adjacent along an edge if and only if their image by $\phi$ are.
This is not always possible (for example if you map a triangle to a square), but in cases $1$ and $3$ it is: both tiling spaces are thus conjugated.
Similarly, cases $4$ and $7$ are conjugated; they admit many tilings referred to in the literature as {\em square-triangle tilings}.

\begin{figure}[hbt]
\centering
\includegraphics[width=\textwidth]{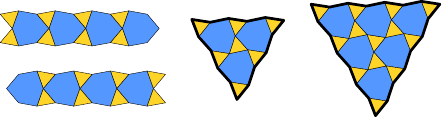}
\caption{
Case $2$ of Fig.~\ref{fig:2}: either two mirrored infinite stripes can be alternated (left) or arbitrarily large ``crested triangles'' can be used to tile as in the triangular grid (center and right).
}
\label{fig:4}
\end{figure}

The case $2$, more original is called {\em shield} in \cite{FHS21} and corresponds to the shield tilings described in \cite{FS24}, see Fig.~\ref{fig:4}.

Now that we have introduced all these notions and illustrated them on the case of triangulated binary packings, we can finally ask our open question:

\begin{question}
What are the conjugacy classes of the triangulated packings with $n$ disk sizes?
How many of them are there?
\end{question}

\begin{figure}[hbt]
\centering
\includegraphics[width=\textwidth]{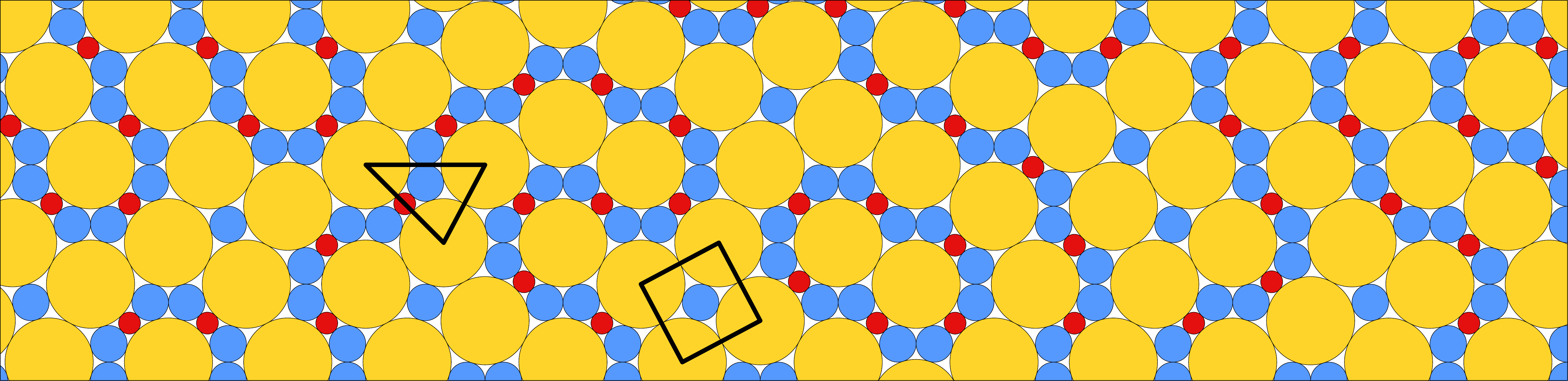}
\caption{
A triangulated ternary packing that can be seen as a tiling by three tiles (the two depicted polygons and the triangle between the centers of three large disks, which does not appear here).
}
\label{fig:5}
\end{figure}

In \cite{FHS21}, the case of three disks is treated.
Surprisingly, all $164$ ternary cases correspond to conjugacy classes that already appear among the $9$ binary cases, except only one shown in Fig.~\ref{fig:5}.
The number of conjugation classes thus seems to grow much slower than the number of cases.

\section{Density}
\label{sec:3}

The {\em density} of a disk packing is defined as the limit superior, when $r$ tends to infinity, of the proportion of the square $[−r, r]^2$ covered by the disks.
This is an important parameter, especially if we think about material science modeling, where the density of a packing of spherical or cylindrical particles plays a key role.
There is never uniqueness of the densest packing (removing a negligible proportion of disks does not change the density) and it can be shown that the density is always reached.

When there is only one disk size, the maximum density is reached for the hexagonal compact packing \cite{FT43, CW10}.
When there are several sizes the maximum density may be higher.
For example, a sufficiently small disk can be inserted in the holes as in the binary packings $8$ and $9$ in Fig.~\ref{fig:1}.
What for larger disks?
For more sizes?
In full generality:

\begin{question}
What does the function that maps $(r_1 , \ldots , r_n )$ onto the maximum density of packings by disks of sizes $r_1 , \ldots , r_n$ look like?
\end{question}

\begin{figure}[hbt]
\centering
\includegraphics[width=\textwidth]{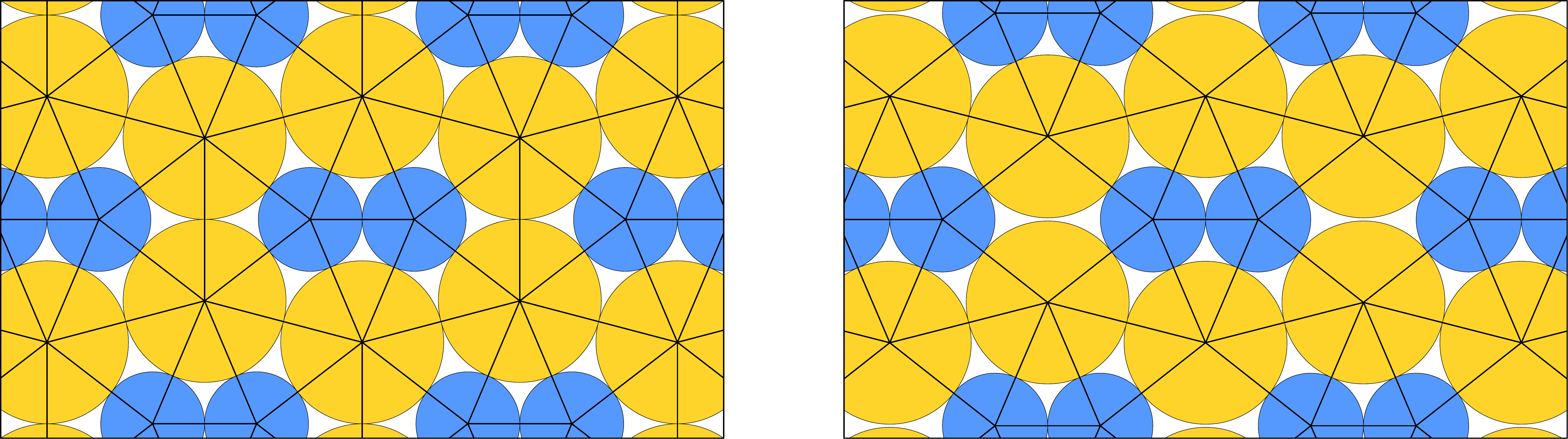}
\caption{
Starting from a triangulated binary packing (left), we release some contacts between disks to allow the small disc to be inflated (right).
Most of the contacts are kept in the hope that the density decreases slowly.
Here, the ratio of disk sizes is increased by about $1.278\%$ and the density decreased to that of the hexagonal compact packing.
}
\label{fig:6}
\end{figure}

\begin{figure}[hbt]
\centering
\includegraphics[width=\textwidth]{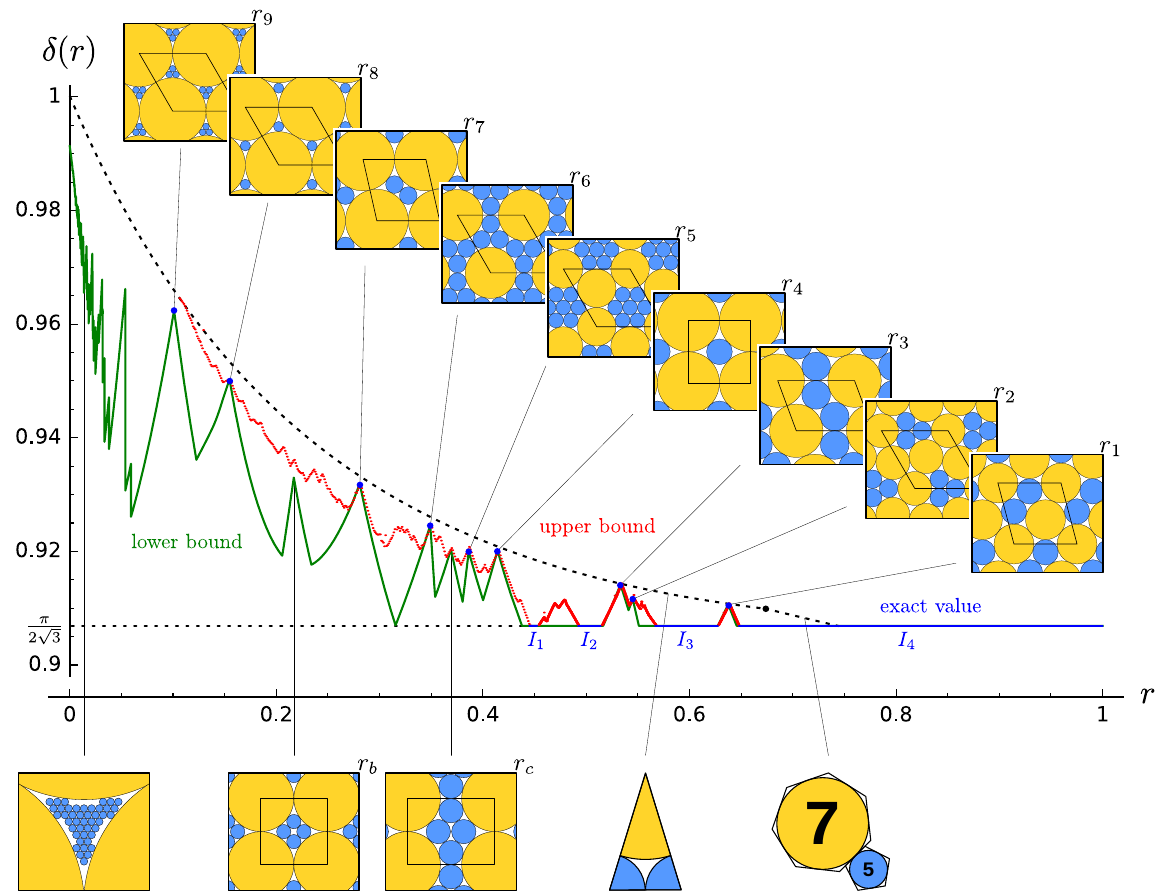}
\caption{
Lower and upper bound on the maximum density $\delta(r)$ of packings of disks of size $1$ and $r < 1$.
The exact value is known for the nine triangulated binary packings as well as for four intervals, where the density of the hexagonal compact packings cannot be improved.
}
\label{fig:7}
\end{figure}

At least, this function is continuous.
Indeed, continuously deflating disks to adjust the ratios, without moving their centers, shows that the maximum density cannot decrease abruptly.
The ratios can actually be adjusted in an even better way, i.e. by decreasing the density less.
The idea, introduced in \cite{FT64} and known as {\em flip and flow} (see, in particular, \cite{CG21}), consists in moving the centers at the same time as the ratios are adjusted, so as to maintain as much contacts as possible between disks, see Fig.~\ref{fig:6}.
A systematic use of this method on particularly dense packings allows to obtain good lower bounds, see for example Fig.~\ref{fig:7} for the case of two disk sizes.

While exhibiting a good packing is sufficient to give a lower bound, obtaining an upper bound is more difficult.
A first result is that of Florian \cite{Flo60}, who bounded from above the maximum density by the density of a triangle connecting the centers of one largest disk and two smallest disks (the density of
a triangle or another finite region is defined as the proportion of its area covered by disks).
A second result is that of Blind \cite{Bli69}, who bounded from above the maximum density by the density in the figure formed by the disjoint union of a regular heptagon circumscribed to a largest disk and a regular pentagon circumscribed to a smallest disk.
Florian’s bound is better below the ratio $r\approx 0.673$, Blind’s beyond.
Both bounds are illustrated for two disk sizes by a dashed line in Fig.~\ref{fig:7}.

Another method to obtain upper bounds is the one called ``localization'' by Lagarias in \cite{Lag11} (part I) and successfully used by Hales and Ferguson to prove the Kepler conjecture (see \cite{Lag11}, parts II–IV).
Basically, given a packing, the idea is to proceed in three steps:
\begin{enumerate}
\item partition the plane or space into compact cells adapted to the packing;
\item consider the density in these cells and ``redistribute'' it between neighboring cells (the denser cells ``give'' density to their less dense neighbors);
\item prove that after redistribution, the density of each cell is lower than the conjectured maximum density.
\end{enumerate}

If this can be done for all packings then the result is proven.
Of course, the difficulty is that there is a continuum of packings which each has infinitely many disks.
To face this, two key ideas are:
\begin{enumerate}
\item restrict the variety of cells by using ``classical'' decompositions (as Voronoi polytopes or Delaunay simplices) in order to reduce the combinatorics of the redistribution process to a finite set of cases;
\item use {\em interval arithmetic} to prove inequalities on compacts (in particular density inequalities on clusters of neighboring cells): the compact is subdivided into a finite number of sufficiently small blocks so that a computation in interval arithmetic on each block allows to conclude.
\end{enumerate}

Following these general principles, each of the nine periodic triangulated binary packings shown in Fig.~\ref{fig:1} has been proven to maximize the density among all the packings with the same disk sizes.
This was first done for Case $4$ by Heppes \cite{Hep00}, who then extended his methods to Cases $1$, $3$, $6$, $7$ and $8$ \cite{Hep03}.
Case $2$ was treated by Kennedy \cite{Ken04} and the two remaining cases in \cite{BF22} which also proposed a unified proof for all the nine cases.

\begin{figure}[hbt]
\centering
\includegraphics[width=0.32\textwidth]{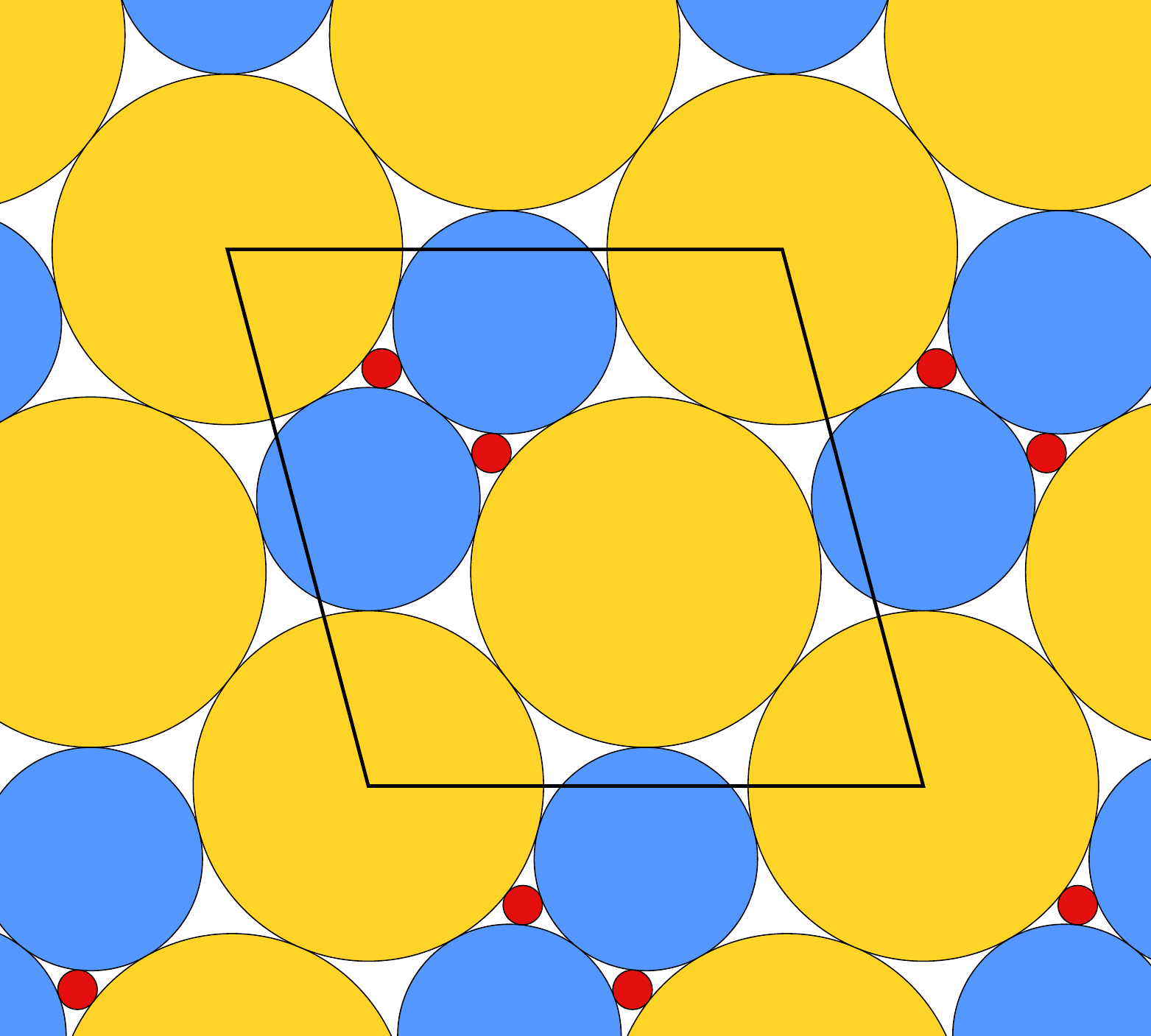}
\hfill
\includegraphics[width=0.32\textwidth]{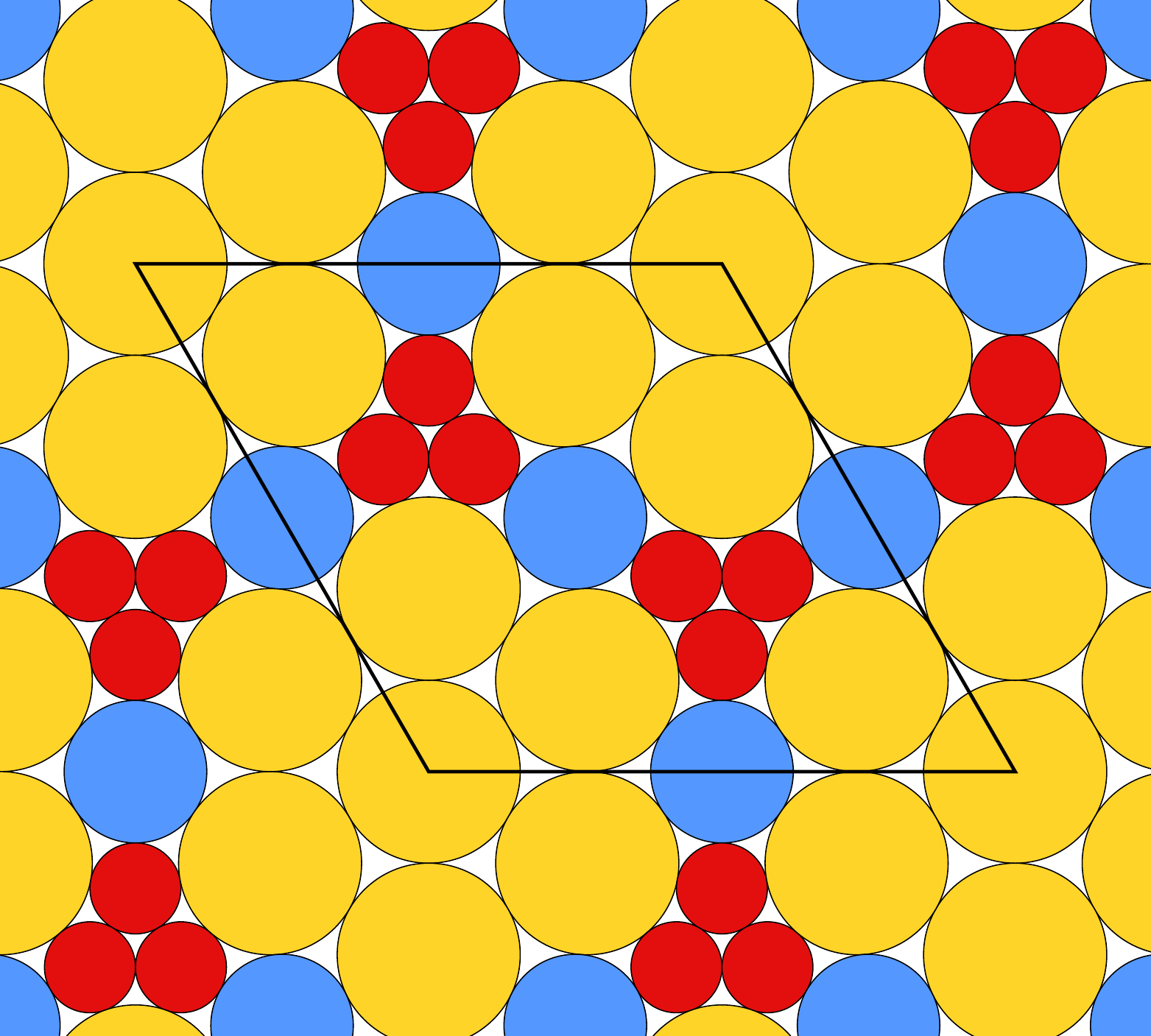}
\hfill
\includegraphics[width=0.32\textwidth]{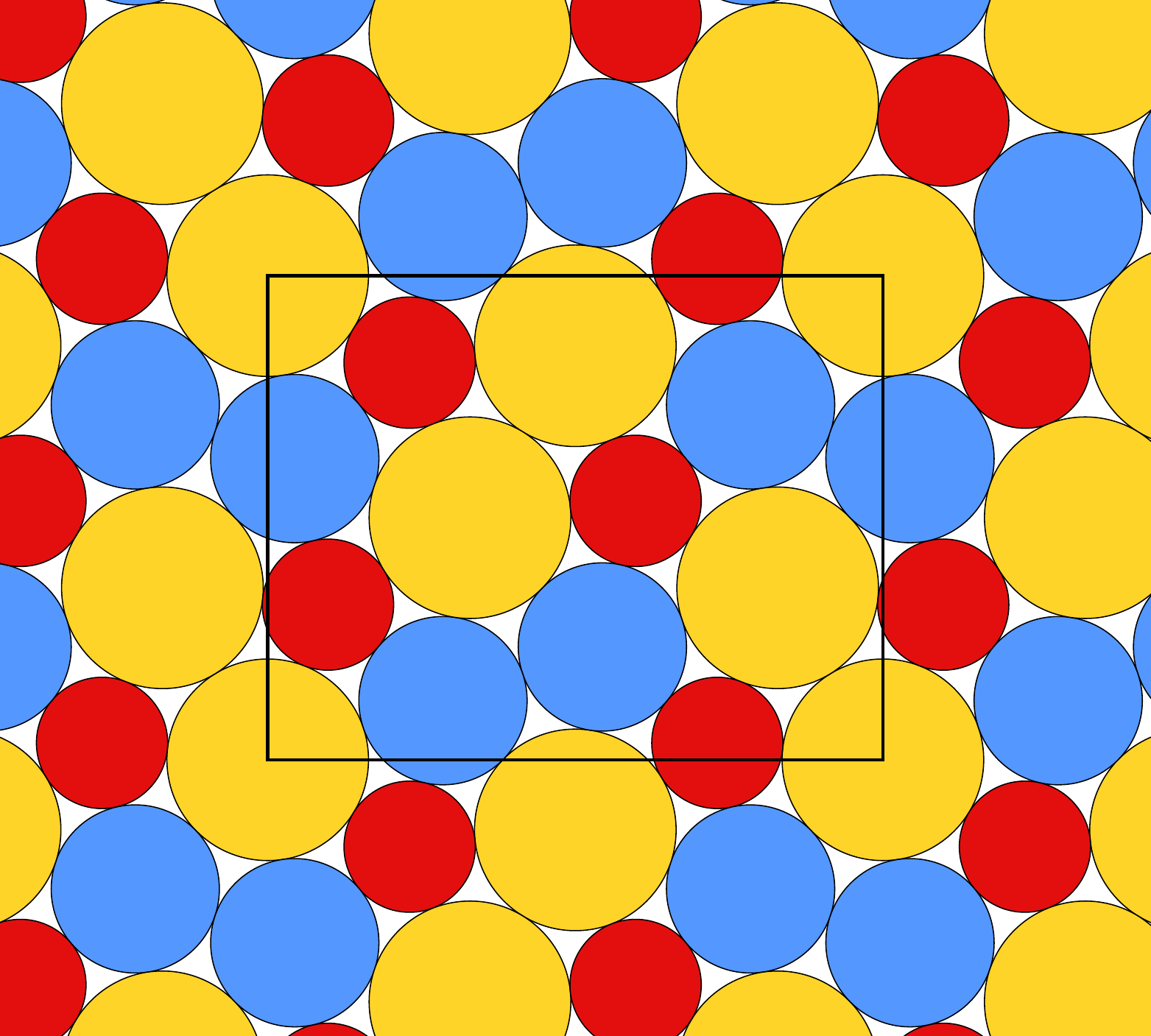}
\caption{Triangulated ternary packings numbered $28$, $110$ and $53$ in \cite{FHS21}.}
\label{fig:8}
\end{figure}

Thus, as soon as two disk sizes allow a triangulated packing, the density among all packings of these disks is maximized over triangulated packings only.
Does this still hold for more disk sizes? 
he densest triangulated packing must be assumed to be {\em saturated}, that is, no disk can be inserted in a hole between three other disks, otherwise there are denser nontriangulated packing (Fig.\ref{fig:8},
left).
Unfortunately, this does not suffice, as illustrated by the packing depicted in Fig.~\ref{fig:8}, center.
In this case, the disks have size $1$, $r\approx 0.779$ and $s \approx 0.497$, where $r$ and $s$ are roots of, respectively:
$$
9x^8 + 12x^7 − 242x^6 + 436x^5 + 665x^4 − 2680x^3 + 2680x^2 − 1056x + 144,
$$
$$
9x^8 − 120x^7 − 380x^6 + 2056x^5 + 12846x^4 − 29672x^3 + 15220x^2 − 2088x + 81.
$$

The point is that the ratio $q := s/r \approx 0.6378$ is very close to the ratio $r_1 \approx 0.6375$ of the binary triangulated packing numbered $1$ in Fig.~\ref{fig:1}.
This latter packing has density greater than $0.9106$ and by slightly inflating its small disks as in Fig.~\ref{fig:6} until the ratio $q$ is reached, the density decreases but remains greater than $0.9105$, whereas the density of the ternary packing depicted in Fig.~\ref{fig:2}, center, is less than $0.9104$.
This example shows that the maximum density is not always reached for a triangulated packing.

The previous example is one of the 40 such cases listed in \cite{FP23} among the $164$ ratios which allow a ternary triangulated packing.
Nevertheless, still in \cite{FP23}, the density has been proven to be maximized over triangulated packings in $16$ cases, including the one in Fig.~\ref{fig:8}, right.
The remaining cases include $15$ nonsaturated cases, $18$ cases where the density is maximized for a binary packing and $75$ open cases.

What about the situation beyond the triangulated packings?
In \cite{Fer22}, the above sketched method of localization has been automated to address all ratios of two disk sizes greater than $0.11$ (for smaller ratios the Florian bound is better).
This is illustrated in Fig.~\ref{fig:7}.

Besides the nine ratios corresponding to the triangulated binary packings, denoted by $r_1 , \ldots , r_9$ in Fig.~\ref{fig:7}, the exact value of the maximum density is known on $4$ intervals.
On these intervals, the densest packing is again triangulated since it is the hexagonal compact packing (for $r \geq  0.74$ it already follows from the Blind bound).

The curve representing the upper bound has many small peaks: each one corresponds to a ratio which locally allows a good arrangement of the disks (typically: an exact corona) but which does not always extend to the whole plane.
At least it does for $r_1 , \ldots , r_9$ , but this does not seem to be the only case.
In particular, the ratios $r_b \approx 0.217$ (root of $x^4 − 4x^3 − 2x^2 − 4x + 1$) and $r_c\approx 0.369$ (root of $x^3 − 5x^2 − x + 1$) allow packings that, without being completely triangulated, have a density that is close to the upper bound (the gap is less than $0.01\%$ for $r_c$).
Similarly, looking around the ratio $r_1 \approx 0.637$ in Fig.~\ref{fig:7} shows that the upper bound is very close to the lower bound, which is obtained by flip and flow (cf Fig.~\ref{fig:6}).
We conjecture:

\begin{conjecture}
For the ratios $r_b \approx 0.217$ and $r_c \approx 0.369$, the periodic (nontriangulated) binary packing depicted in Fig.~\ref{fig:7} maximizes the density among all the binary packings of disks within the same size ratio.
Further, for every ratio around $r_1 \approx 0.637$ such that the relaxation of the triangulated binary packing depicted in Fig.~\ref{fig:6} is denser than the hexagonal compact packing (this includes the interval $[0.627, 0.645]$), then this relaxation maximizes the density.
\end{conjecture}

It would also be interesting to obtain an equivalent of Fig.~\ref{fig:7} for three disk sizes (the maximum density is then a surface in $\mathbb{R}^3$).
Lower bounds could be obtained by efficiently applying the flip and flow method to triangulated ternary packings of \cite{FHS21}, while upper bounds could be obtained by automating the techniques of \cite{FP23} for all possible size ratios.

\section{Uniformity}

\begin{figure}[hbt]
\centering
\includegraphics[width=\textwidth]{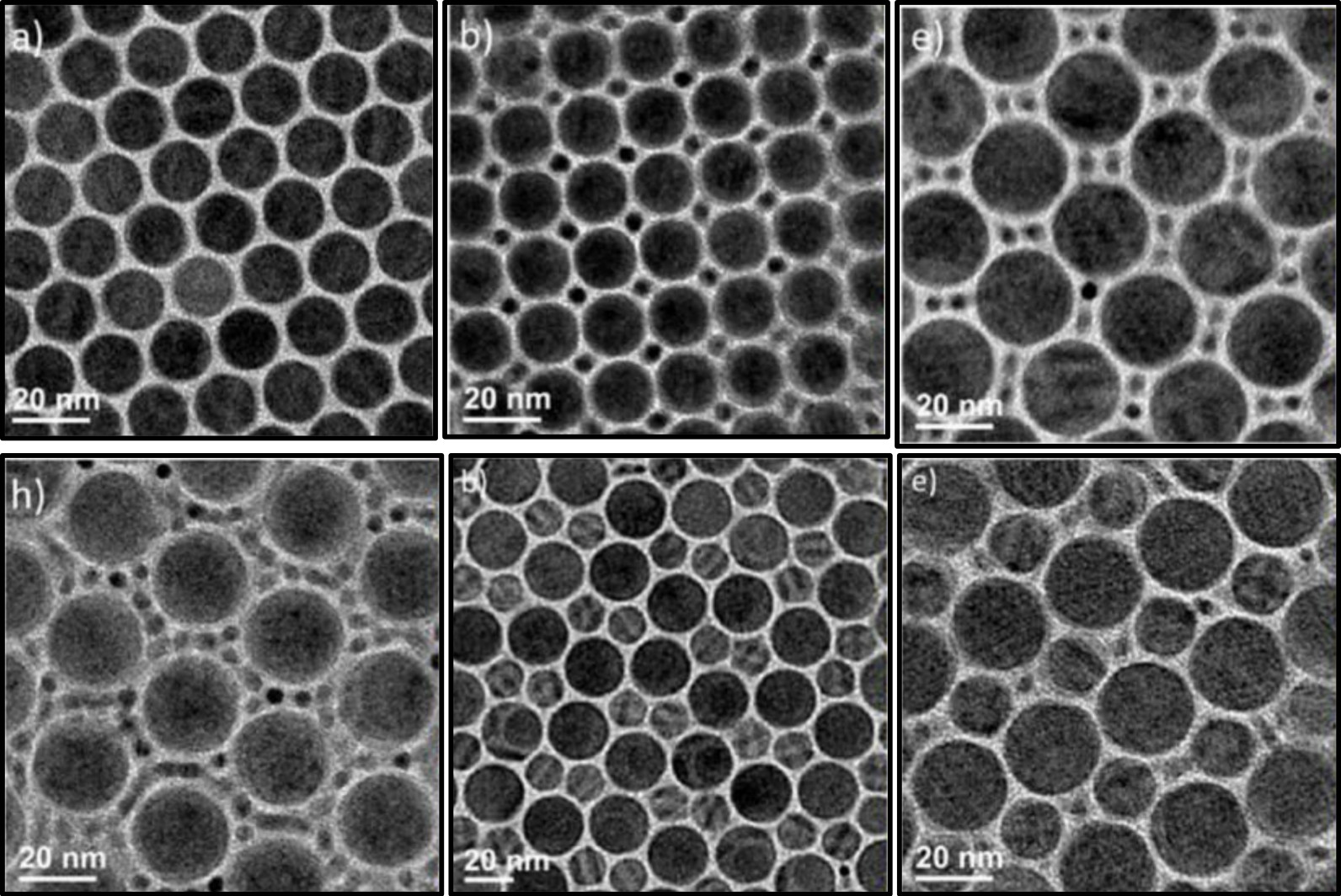}
\caption{
Self-assembled nanorods reported in \cite{PDKM15} (courtesy of the authors).
In the reading direction, we can recognize the hexagonal compact packing, the triangulated binary packings numbered $4$, $7$, $6$ and $1$ in Fig.~\ref{fig:1} as well as a triangulated ternary packing (numbered $82$ in \cite{FHS21}.
}
\label{fig:9}
\end{figure}

A packing can always be densified by inserting small discs in the holes.
What if, conversely, the disks are forced to be of ``comparable'' size?
This question is of particular interest for applications in materials science, where we want to combine particles of the same order of magnitude.
The resulting structure must be denser than a phase separation of each particle size, with the idea being that the maximization of the density is favored by physical constraints (such as attractive forces between particles).
For example, \cite{PDKM15} reports the experimental synthesis of self-assembled nanorods which yield structure looking exactly the same as triangulated disk packings, see Fig.~\ref{fig:9}.

Formally, the {\em uniformity} of a packing of disk is defined as the ratio of the smallest to the largest disk sizes.
Since any uniformity can be achieved with two sizes of disks by combining hexagonal compact packings (HCP) of each size on either side of a line, we restrict to packings denser than HCP and ask:

\begin{question}
What is the highest uniformity among the packings denser than the HCP?
What if we consider only triangulated packings or those that maximize density for given disk sizes?
\end{question}

For a long time, the record holders where the packings depicted in Fig.~\ref{fig:6}, respectively in the ``triangulated'' category on the left (uniformity $u \approx 0.6376$) and in the ``free'' category on the right (uniformity $u \approx 0.6457$) \cite{FT64}.
Soon after the full list of triangulated ternary packings was established, Connelly spotted the one numbered $53$, depicted in Fig.~\ref{fig:8}, right, which improved the maximal uniformity to $u \approx 0.6510$ root of
$$
89x^8 + 1344x^7 + 4008x^6 − 464x^5 − 2410x^4 + 176x^3 + 296x^2 − 96x + 1.
$$
With his student, they deformed it by flipping and flowing to reach uniformity $u \approx 0.6585$ at the cost of the triangulated hypothesis \cite{CP}.
The deformation is hardly visible, let us just say that it is obtained by releasing the contact between the larger disks and inflating the smaller disks, while the other two sizes remain fixed.
Connelly subsequently proposed the following conjecture:

\begin{conjecture}
The ternary triangulated packing numbered $53$ in \cite{FHS21} and depicted in Fig.~\ref{fig:8}, right, is the most uniform triangulated packing (HCP excepted).
Its deformation proposed in \cite{CP} yields the highest possible uniformity for a disk packing (HCP excepted).
\end{conjecture}

\begin{figure}[hbt]
\centering
\includegraphics[width=0.24\textwidth]{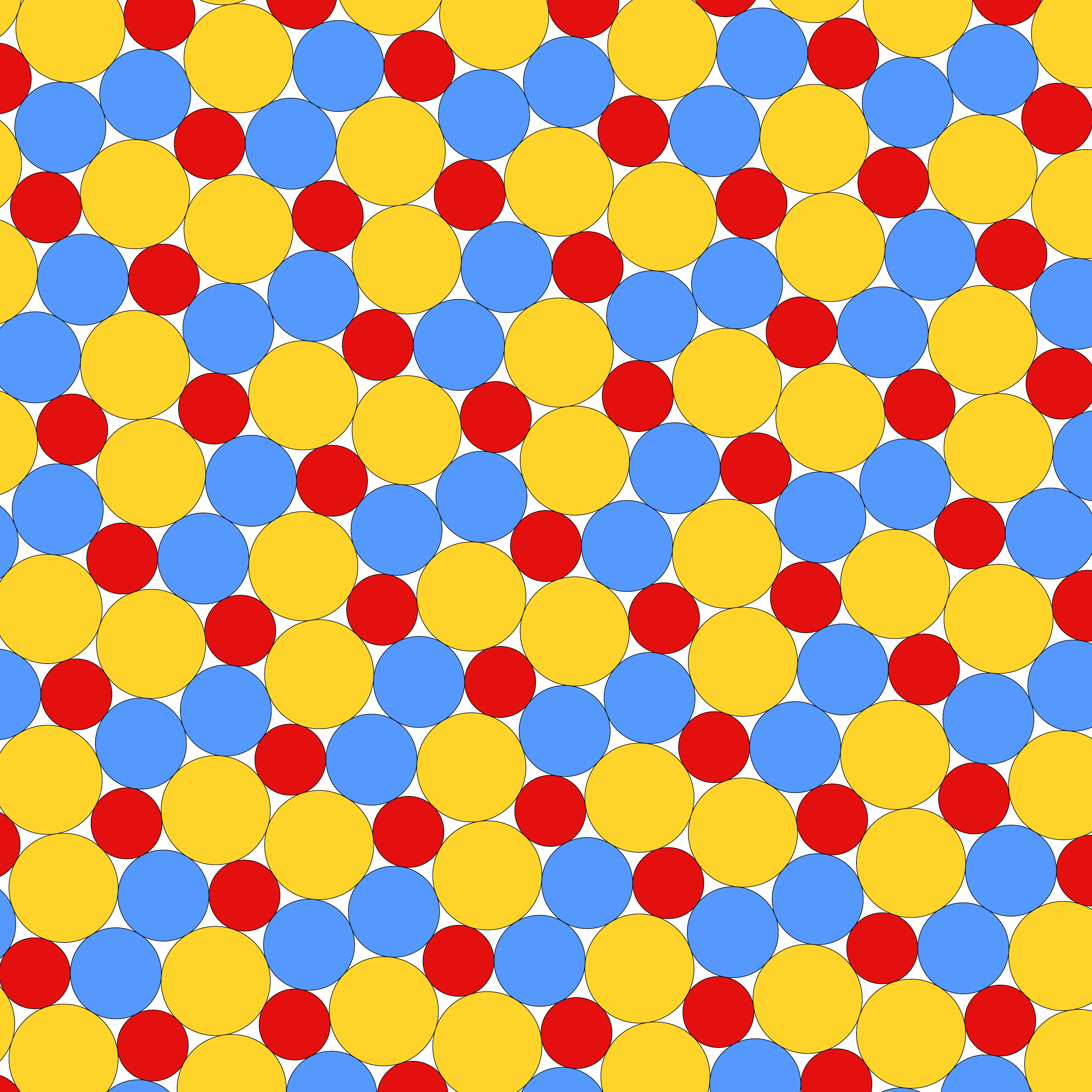}
\hfill
\includegraphics[width=0.24\textwidth]{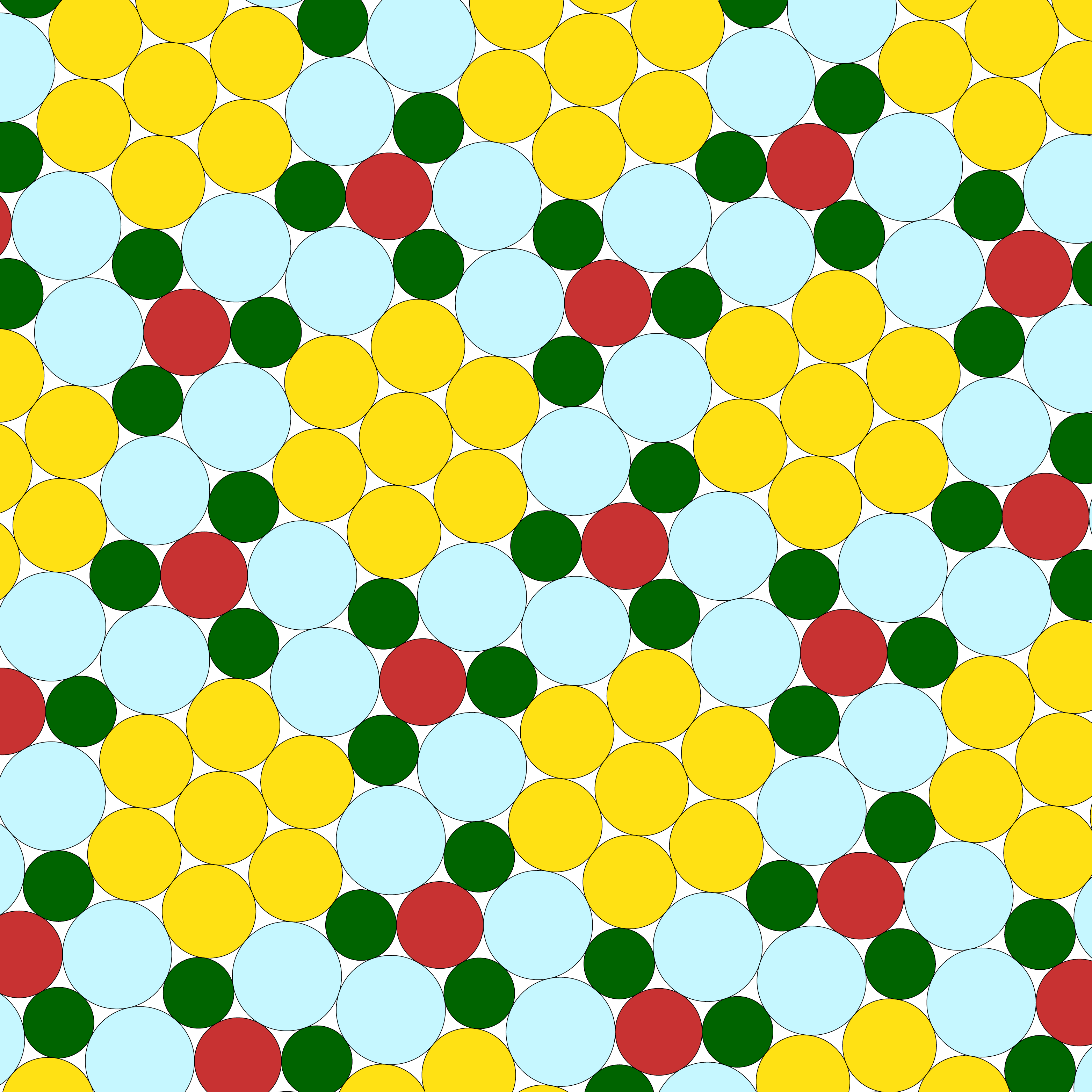}
\hfill
\includegraphics[width=0.24\textwidth]{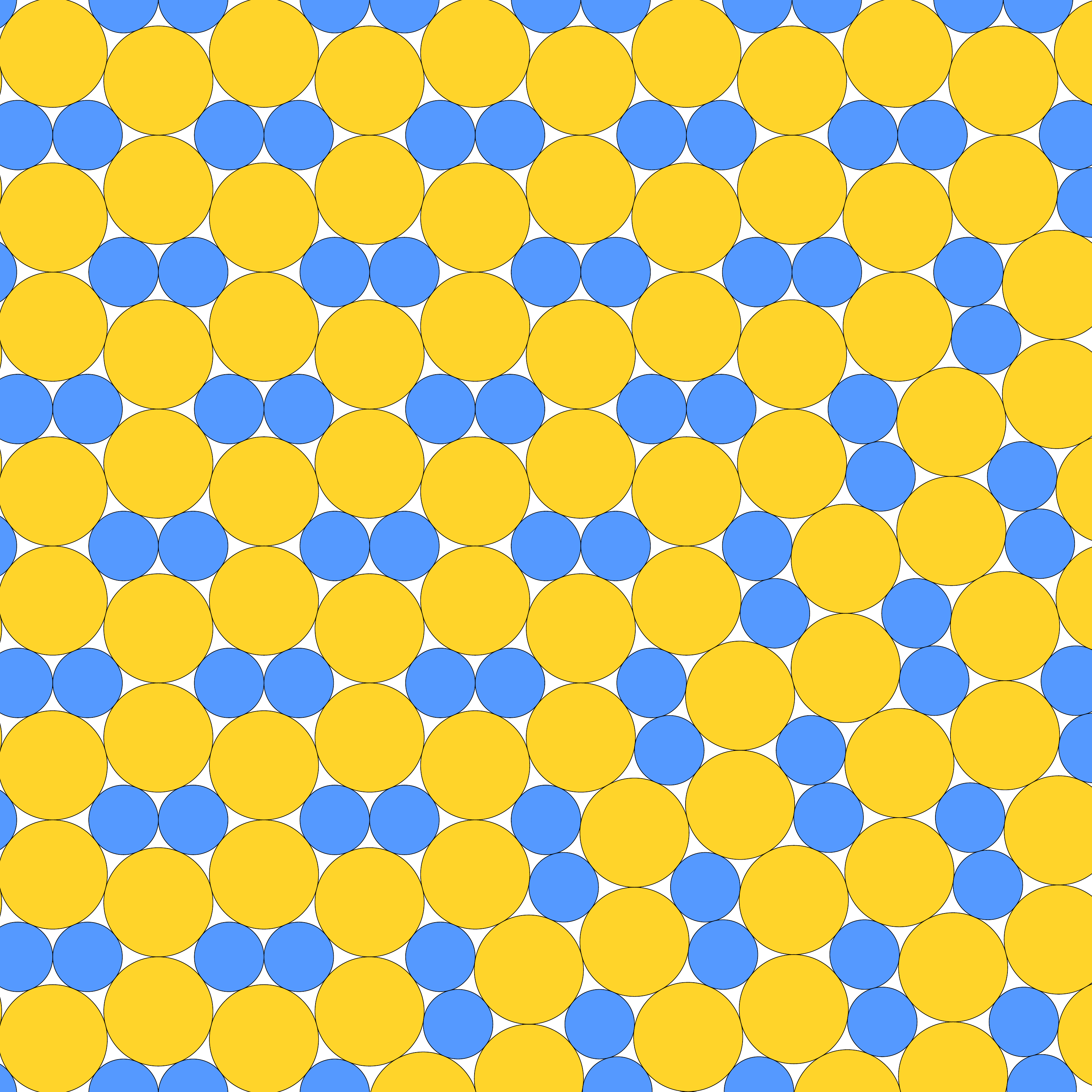}
\hfill
\includegraphics[width=0.24\textwidth]{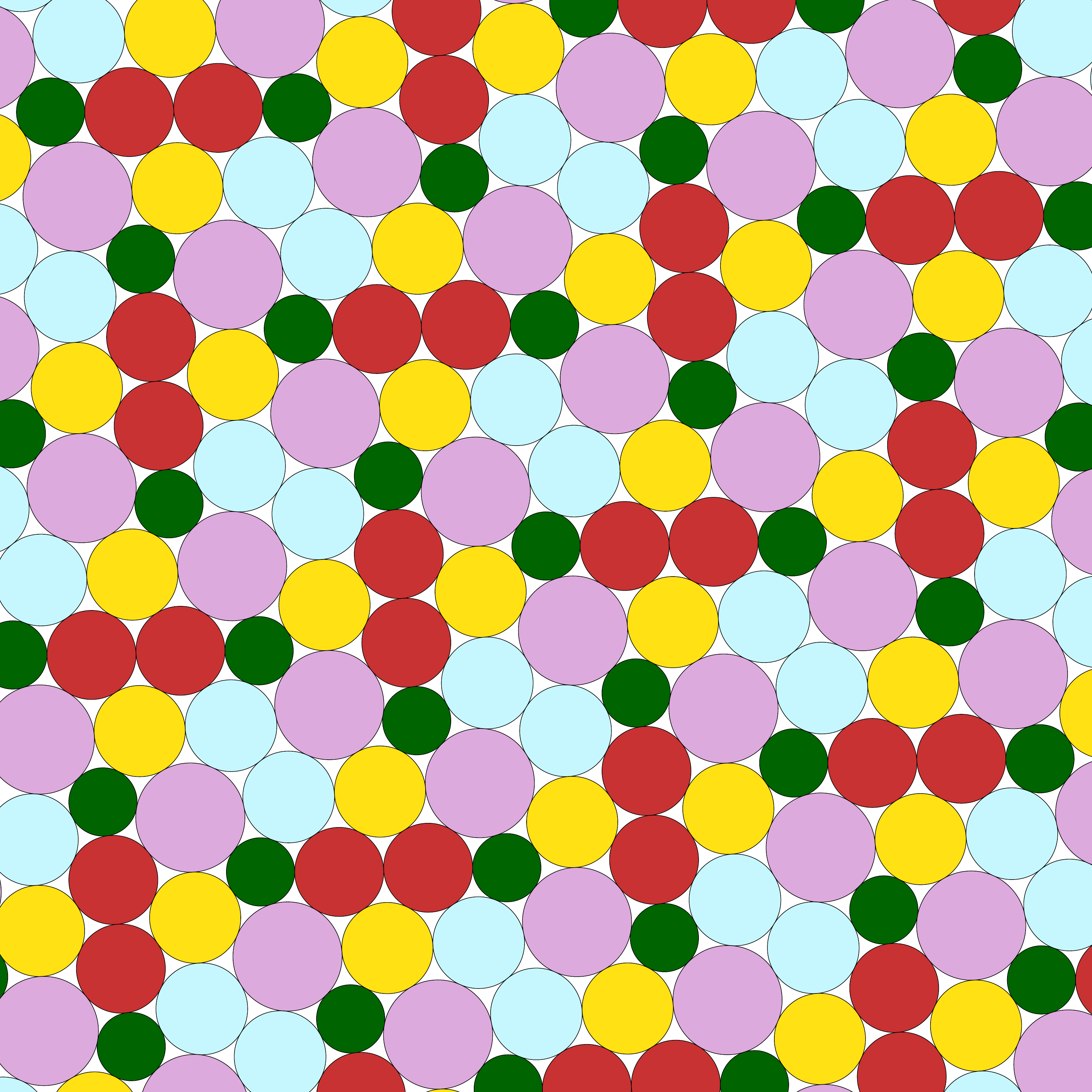}
\caption{
The most uniform triangulated packings with at most $5$ disk sizes, in decreasing order of uniformity from left to right.
}
\label{fig:10}
\end{figure}

This may seem doubtful because considering more disk sizes should probably further improve uniformity.
However, there might be good reasons to think that it holds.
For example, an exhaustive search for triangulated packings with uniformity at least $0.6$ (this assumption reduces a lot the computation time, which however increases quickly with the number of disks) shows that uniformity is less than $0.6486$ for $4$ disk sizes and less than $0.6258$ for $5$ disk sizes (Fig.~\ref{fig:10}).
It is unlikely that deformations on these packings will do better than the one proposed in \cite{CP}.
An interesting approach based on the notion of rigidity of a structure has been developed in \cite{CZ}, without however being able to prove this conjecture to date.

Conversely, the Blind bound mentionned in the previous section yields (by finding the ratio such that the density in the union of the heptagon and the pentagon is that of HCP) that the uniformity of any packing denser than HCP is at most
$$
\sqrt{\frac{7\tan\tfrac{\pi}{7}-6\tan\tfrac{\pi}{6}}{6\tan\tfrac{\pi}{6}-5\tan\tfrac{\pi}{5}}}\approx 0.74299.
$$
In the case of triangulated packings, it is much easier to get a better upper bound.
Indeed, the smallest disk has at most $5$ neighbors (otherwise the packing is HCP), and its size is maximized for the largest possible neighbors.
The maximum uniformity is thus bounded from above by the size ratio of a disk surrounded by five equally sized disks.
A computation yields the bound
$$
\frac{2\sqrt{2}}{\sqrt{5-\sqrt{5}}}-1\approx 0.70130.
$$
There is still a gap between this upper bound and the lower bound of approximatively $0.6510$ given by the best ternary triangulated packing.

\section{Complexity}
\label{sec:5}

In Section~\ref{sec:1}, we mentionned that the problem of finding all the triangulated packings could become intractable, because given disk sizes which allow to form coronas around disks of each size, it is not trivial to determine whether this can be used to build a triangulated packing over the whole plane.
The underlying problem is actually referred to as the {\em Domino problem}.

\begin{figure}[hbt]
\centering
\includegraphics[width=0.8\textwidth]{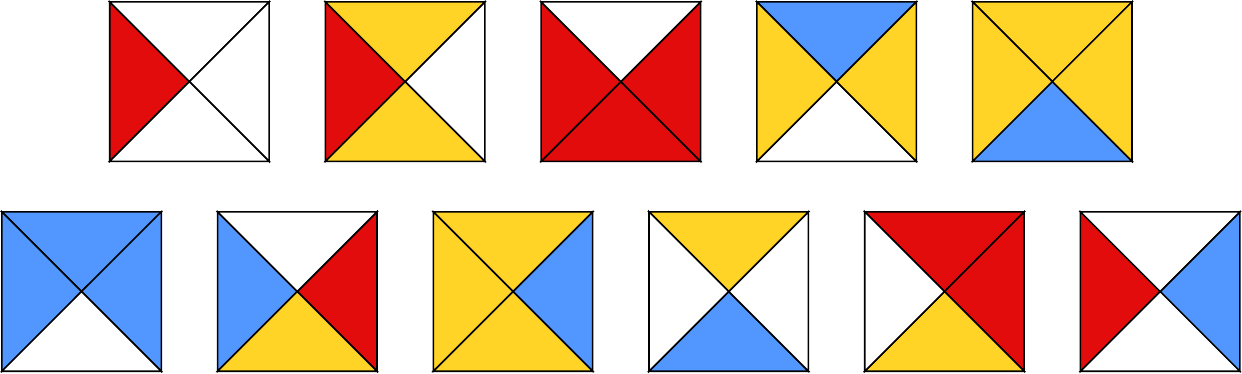}
\caption{These $11$ Wang tiles do admit tiling of the plane, but only non-periodic ones; it is moreover of minimal cardinality with this property \cite{JR15}.}
\label{fig:11}
\end{figure}

The Domino problem can be simply stated using the terminology introduced at the beginning of Section~\ref{sec:2}: given a finite number of unit squares with colored edges, called {\em Wang tiles} and which can only be translated (Fig.~\ref{fig:11}), does there exist a tiling?
This problems has been proved to be undecidable in the 1960’s \cite{Ber66}, that is, there is no algorithm which decides the existence of a tiling in finite time for any given set of Wang tiles.
Very roughly, the idea is to show that any Turing machine can be simulated by a finite set of Wang tiles: the tape at time $t$ correspond to the $t$-th line of the tiling, and there exists a tiling of the whole
plane if and only if the Turing maching does not halt.
The undecidability of the Domino problem then follows from the undecidability of the {\em Halting problem}, proved by Turing in 1937.

The above proof sketch, however, omits a key argument: the existence of an {\em aperiodic tile set}, that is, a finite set of tiles which allow tilings of the plane but only non-periodic ones (no invariance by translation).
Indeed, if any set of tiles that tile the plane could do this periodically, the Domino problem could be solved as follows: try all the ways to form a $n \times n$ square for larger and larger $n$ (this requires a huge but finite amount of time).
If there is no tiling, then for $n$ large enough there will be no such square and the answer will be negative.
If there is a periodic tiling, then for $n$ large enough there will be a square whose left and right sides (resp. top and bottom) are identical and which thus allows to periodically tile the plane: the answer will be positive.

The existence of an aperiodic tile set is used in a subtle way in the proof of the Domino problem.
Basically, a suitable aperiodic set of tiles is combined with the tiles that simulate the Turing machine so as to force the computation to start somewhere (without the same tape would repeat forever on each line
of the tiling).
It is a little more complicated in the details, because the same computation is then necessarily started everywhere in the tiling (otherwise there are arbitrarily large areas without computation, so by compactness a tiling of the whole plane without any computation).
An elegant and comprehensive account can be found in \cite{Rob71} (it relies on the so-called {\em Robinson aperiodic tile set}).

What does this have to do with disk packings?
The problem of the existence of a triangulated disk packing, addressed in Section~\ref{sec:1}, turns out to be a variant of the Domino problem with a very particular class of tiles: the triangles of a
triangulated packing.
The question then arises: does this change the nature (decidable or not) of the problem? Formally:

\begin{question}
Does exist an algorithm which decides in finite time whether a given finite set of disk sizes allows a triangulated packing other than the hexagonal compact packing?
\end{question}

The disk sizes must be assumed to be at least computable so that they can be considered as an input to the algorithm.
In fact they can even be assumed to be algebraic because it is a necessary condition for the existence of a triangulated packing.
Each size can thus be given by an integer polynomial of which it is a root, with a rational interval that contains this root.

The above question is of course interesting in connection with the search for triangulated packings, but it also has its own interest in computability theory.
Indeed, it is not very clear how aperiodicity is the key to the undecidability of the Domino problem: it could be that the aperiodic tile set used in the proof is very particular and that for a more restricted class of tile sets, the problem would become decidable.
Jeandel proved in \cite{Jea10} that this was not the case for a related problem, the {\em Periodic Domino problem}, which asks whether a set of tiles admits a periodic tiling: its proof is based on the mere existence of an aperiodic tile set, without any particular hypothesis about it.
Anyway, the first step is to find ``aperiodic disks'', otherwise the problem can be decided like the general variant with Wang tiles:

\begin{conjecture}
There exists a finite set of disk sizes which allow triangulated packings other than the hexagonal compact packing, with no such packing being periodic (weaker version: with the densest one being not periodic).
\end{conjecture}

What strategy should be considered to find such disc sizes?
A ``bottom-up approach'' would be to classify all triangulated packings of $k$ disk sizes for increasing values of $k$, as discussed in Sections~\ref{sec:1} and ~\ref{sec:2}.
This is the approach successfully followed in \cite{JR15}, which led to the smallest set of aperiodic Wang tiles (Fig.~\ref{fig:11}), though the problem is here made more complicated by the geometry of the tiles.
Alternatively, a ``top-down approach'' would be to start from a known non-periodic tiling and to derive a disk packing for which this tiling would be the contact graph.
To this end, one could consider starting from a hexagonal compact packing then flip and flow as in \cite{CG21} until transforming the contact graph into the desired non periodic graph.

\bibliographystyle{alpha}
\bibliography{survey}

\end{document}